\documentclass[12pt]{article}
\usepackage{epic,eepic}
\usepackage[dvips]{epsfig}
\usepackage{color}
\usepackage{amsmath, amscd, amssymb}
\usepackage{amssymb}
\def\draw #1 by #2 (#3){
  \vbox to #2{
    \hrule width #1 height 0pt depth 0pt
    \vfill
    \special{picture #3}
    }
  }
\def\scaleddraw #1 by #2 (#3 scaled #4){{
  \dimen0=#1 \dimen1=#2
  \divide\dimen0 by 1000 \multiply\dimen0 by #4
  \divide\dimen1 by 1000 \multiply\dimen1 by #4
  \draw \dimen0 by \dimen1 (#3 scaled #4)}
  }

\def \vski{\vspace{12pt}}
\newcommand{\ind}{\hspace{.25in}}
\def \bea{\begin{eqnarray}}
\def \eea{\end{eqnarray}}

\begin{document}
\renewcommand{\labelenumi}{\theenumi}
\newcommand{\qed}{\mbox{\raisebox{0.7ex}{\fbox{}}}}
\newtheorem{theorem}{Theorem}
\newtheorem{example}{Example}
\newtheorem{problem}[theorem]{Problem}
\newtheorem{defin}[theorem]{Definition}
\newtheorem{lemma}[theorem]{Lemma}
\newtheorem{corollary}[theorem]{Corollary}
\newtheorem{nt}{Note}
\newtheorem{proposition}[theorem]{Proposition}
\renewcommand{\thent}{}
\newenvironment{pf}{\medskip\noindent{\textbf{Proof}:  \hspace*{-.4cm}}\enspace}{\hfill \qed \medskip \newline}
\newenvironment{pft1}{\medskip\noindent{\textbf{Proof of Theorem 1}:  \hspace*{-.4cm}}\enspace}{\hfill \qed \medskip \newline}
\newenvironment{pft2}{\medskip\noindent{\textbf{Proof of Theorem 2}:  \hspace*{-.4cm}}\enspace}{\hfill \qed \medskip \newline}
\newenvironment{defn}{\begin{defin}\em}{\end{defin}}{\vspace{-0.5cm}}
\newenvironment{lem}{\begin{lemma}\em}{\end{lemma}}{\vspace{-0.5cm}}
\newenvironment{cor}{\begin{corollary}\em}{\end{corollary}}{\vspace{-0.5cm}}
\newenvironment{thm}{\begin{theorem} \em}{\end{theorem}}{\vspace{-0.5cm}}
\newenvironment{pbm}{\begin{problem} \em}{\end{problem}}{\vspace{-0.5cm}}
\newenvironment{note}{\begin{nt} \em}{\end{nt}}{\vspace{-0.5cm}}
\newenvironment{exa}{\begin{example} \em}{\end{example}}{\vspace{-0.5cm}}
\newenvironment{pro}{\begin{proposition} \em}{\end{proposition}}{\vspace{-0.5cm}}

\setlength{\unitlength}{12pt}
\newcommand{\comb}[2]{\mbox{$\left(\!\!\begin{array}{c}
            {#1} \\[-0.5ex] {#2} \end{array}\!\!\right)$}}
\renewcommand{\labelenumi}{(\theenumi)}
\renewcommand{\b}{\beta}
\newcounter{myfig}
\newcounter{mytab}
\def\mod{\hbox{\rm mod }}
\def\scaleddraw #1 by #2 (#3 scaled #4){{
  \dimen0=#1 \dimen1=#2
  \divide\dimen0 by 1000 \multiply\dimen0 by #4
  \divide\dimen1 by 1000 \multiply\dimen1 by #4
  \draw \dimen0 by \dimen1 (#3 scaled #4)}
  }
\newcommand{\Aut}{\mbox{\rm Aut}}
\newcommand{\w}{\omega}
\def\r{\rho}
\newcommand{\DbF}{D \times^{\phi} F}
\newcommand{\autF}{{\tiny\Aut{\scriptscriptstyle(\!F\!)}}}
\def\Cay{\mbox{\rm Cay}}
\def\a{\alpha}
\newcommand{\C}[1]{\mathcal #1}
\newcommand{\B}[1]{\mathbb #1}
\newcommand{\F}[1]{\mathfrak #1}
\title{\textbf{There are only finitely many distance-regular graphs with valency $\textbf{k}$ at least three, fixed ratio $\mathbf{\frac{k_2}{k}}$ and large diameter}}
  \author{Jongyook Park$^{\,\rm 1}$
,~  Jack H. Koolen$^{\,\rm 1}$
and Greg Markowsky$^{\,\rm 2}$\\
{\small {\tt jongyook@postech.ac.kr} ~~
{\tt koolen@postech.ac.kr} ~~ {\tt gmarkowsky@gmail.com}}\\
{\footnotesize{$^{\rm 1}$Department of Mathematics,  POSTECH, Pohang 790-785, South Korea}}\\
{\footnotesize{$^{\rm 2}$Pohang Mathematics Institute,  POSTECH, Pohang 790-785, South Korea}}}
\maketitle

\begin{abstract}

In this paper, we show that for given positive integer $C$, there are only finitely many distance-regular graphs with valency $k$ at least three, diameter $D$ at least six and $\frac{k_2}{k}\leq C$. This extends a conjecture of Bannai and Ito.

\\
\end{abstract}

\section{Introduction}

In 1984, Bannai and Ito made the following conjecture (Cf. \cite[p.237]{bi}).

\vski

\textbf{Conjecture:} There are finitely many distance-regular graphs with fixed valency at least three.

\vski

This conjecture has been proved in \cite{bdkm}. In this paper we strengthen this result as follows.

\begin{thm} \label{2}
Let $C$ be a positive integer. Then there are finitely many distance-regular graphs with valency $k>2$, diameter $D\geq6$ and $\frac{k_2}{k}\leq C$.
\end{thm}

\textbf{Remark:} The Hadamard graphs of order $2\mu$ (\cite[Section 1.8]{bcn}) have intersection array $\{2\mu,2\mu-1,\mu,1;1,\mu,2\mu-1,2\mu\}$ and have $\frac{k_2}{k}=\frac{2\mu-1}{\mu}<2$. For each Hadamard matrix of order $2\mu$, there exists a Hadamard graph of order $2\mu$, and for each $\mu=2^n$ for $n\geq0$, there exists a Hadamard matrix of order $2\mu$. The Taylor graphs, i.e., antipodal 2-covers of diameter three, have $k_2=k$ and there are infinitely many Taylor graphs (see \cite[Section 1.5]{bcn}). The complement of a non-complete strongly regular graph is also a strongly regular graph. Thus, Theorem~\ref{2} is not true for $D\leq4$. It is not known whether it is true for $D=5$.

\ind We remark further that in \cite{kp} it was shown that if $\frac{k_2}{k}\leq\frac{3}{2}$, then either (D=3 and the graph $\Gamma$ is bipartite, a Taylor graph or the Johnson graph $J(7,3)$) or (D=4 and the graph $\Gamma$ is the 4-cube). Theorem~\ref{2} is an extension of this theorem as well.

\vski

In order to prove Theorem~\ref{2}, we will show

\begin{thm}\label{6}
Let $C$ be a positive integer. Then there are only finitely many distance-regular graphs with valency $k>2$, diameter $D\geq 6$ and $\frac{b_2}{c_2}\leq C$.
\end{thm}

Theorem~\ref{2} immediately follows from this theorem as $\frac{k_2}{k} = \frac{b_1}{c_2}$ and $b_2\leq b_1$. The paper is organized as follows. In Section 2 we will give the definitions and basic facts required for the ensuing arguments. In Section 3 we give a lower bound on the second largest eigenvalue, and in Section 4 we will discuss distance-regular Terwilliger graphs. In Section 5 we consider distance-regular graphs with a fixed $\frac{b_t}{c_t}$ for some $t\geq1$ and give a proof of Theorem~\ref{6}.

\section{Definitions and preliminaries}

All the graphs considered in this paper are finite, undirected and
simple (for unexplained terminology and more details, see \cite{bcn}). Suppose
that $\Gamma$ is a connected graph with vertex set $V(\Gamma)$ and edge set $E(\Gamma)$, where $E(\Gamma)$ consists of unordered pairs of two adjacent vertices. The {\em distance} $d(x,y)$ between
any two vertices $x$ and $y$ of $\Gamma$
is the length of a shortest path connecting $x$ and $y$ in $\Gamma$. We denote $v$ as the number of vertices of $\Gamma$ and define the {\em diameter} $D$
of $\Gamma$ as the maximum distance in $\Gamma$.  For a vertex $x \in V(\Gamma)$, define $\Gamma_i(x)$ to be the set of
vertices which are at distance precisely $i$ from $x~(0\le i\le D)$. In addition, define $\Gamma_{-1}(x) = \Gamma_{D+1}(x)
:= \emptyset$. We write $\Gamma(x)$ instead of $ \Gamma_1(x)$ and denote $x\sim_{\Gamma} y$ or simply $x\sim y$ if two vertices $x$ and $y$ are adjacent in $\Gamma$. The {\em adjacency matrix} $A$ of the graph $\Gamma$ is the (0,1)-matrix whose rows and columns are indexed by the vertex set $V(\Gamma)$ and the $(x,y)$-entry is $1$ whenever $x\sim y$ and 0 otherwise.
The {\em eigenvalues} of the graph $\Gamma$ are the eigenvalues of  $A$.\\

 For a connected graph $\Gamma$, the {\em local graph} $\Delta(x)$ at a vertex $x\in V(\Gamma)$ is the subgraph induced on $\Gamma(x)$. Let $\Delta$ be a graph.  If the local graph $\Delta(x)$ is isomorphic to $\Delta$ for any vertex $x\in\Gamma(x)$, then we say $\Gamma$ is {\em locally} $\Delta$.\\

For a graph $\Gamma$, a partition $\Pi= \{P_1, P_2, \ldots ,
P_{\ell}\}$ of the vertex set $V(\Gamma)$ is called {\em
equitable} if there are constants $\beta_{ij}$ such that each vertex $x \in P_i$ has exactly $\beta_{ij}$ neighbors in $P_j$ ($1\leq i, j \leq \ell$). The {\em quotient
matrix} $Q(\Pi)$ associated with the equitable partition $\Pi$ is
the $\ell \times \ell$ matrix whose $(i,j)$-entry $Q(\Pi)_{(i,j)}$ is $\beta_{ij}$ ($1\leq i,j\leq \ell$). Note that the eigenvalues of the quotient matrix $Q(\Pi)$ are also eigenvalues (of the adjacency matrix $A$) of $\Gamma$ \cite[Theorem 9.3.3]{godsil-combin}.\\

A connected graph $\Gamma$ with diameter $D$ is called {\em{distance-regular}} if there are integers $b_i, c_i$ $(0\leq i\leq D)$ such that for any two vertices $x, y \in V(\Gamma)$ with $d(x, y)=i$, there are precisely $c_i$ neighbors of $y$ in $\Gamma_{i-1}(x)$ and $b_i$ neighbors of $y$ in $\Gamma_{i+1}(x)$, where we define $b_D=c_0=0$. In particular, any distance-regular graph  is regular with valency $k := b_0$. Note that a  (non-complete) connected {\em strongly regular graph} is just a distance-regular graph with diameter two. We define $a_i := k-b_i-c_i$ for notational convenience.  Note that $a_i=\mid\Gamma(y)\cap\Gamma_i(x)\mid$ holds for any two vertices $x, y$ with $d(x, y)=i$ $(0\leq i\leq D).$\\

For a distance-regular graph $\Gamma$ and a vertex $x\in V(\Gamma)$, we denote $k_i:=|\Gamma_i(x)|$, and it is easy to see that $k_i = \frac{b_0 b_1 \cdots b_{i-1}}{c_1 c_2 \cdots c_i}$ and hence does not depend on $x$.  The numbers $a_i$, $b_{i-1}$ and $c_i$ $(1\leq i\leq D)$ are called the {\em{intersection~numbers}}, and the array $\{b_0,b_1,\cdots,b_{D-1};c_1,c_2,\cdots,c_D\}$ is called the {\em{intersection~array}} of $\Gamma$. We define

\begin{center}
$h=h(\Gamma):=|\{i~|~1\leq i \leq D-1~~ {\rm and}~~ (c_i,a_i,b_i)=(c_1,a_1,b_1)\}|$.
\end{center}

Some standard properties of the intersection numbers are collected in the following lemma.

\begin{lem}({\rm \cite[Proposition 4.1.6]{bcn}})\label{pre}{\ \\}
Let $\Gamma$ be a distance-regular graph with valency $k$ and diameter $D$. Then the
following holds:\\
$(i)$ $k=b_0> b_1\geq \cdots \geq b_{D-1}~;$\\
$(ii)$ $1=c_1\leq c_2\leq \cdots \leq c_{D}~;$\\
$(iii)$ $b_i\ge c_j$ \mbox{ if }$i+j\le D~.$
\end{lem}

Suppose that $\Gamma$ is a distance-regular graph with valency $k\ge 2$ and diameter $D\ge 1$. Then $\Gamma$ has exactly $D+1$ distinct eigenvalues, namely $k=\theta_0>\theta_1>\cdots>\theta_D$ (\cite[p.128]{bcn}), and  the {\em multiplicity} of
$\theta_i$ ($0\le i\le D$) is denote by  $m_i$. For an eigenvalue $\theta$ of $\Gamma$, the sequence $(u_i)_{i=0,1,...,D}$ = $(u_i(\theta))_{i=0,1,...,D}$
satisfying $u_0$ = $u_0(\theta)$ = $1$, $u_1$ = $u_1(\theta)$ = $\theta/k$, and
\begin{center}
$c_i u_{i-1} + a_i u_i + b_i u_{i+1} = \theta u_i$ $(i=1,2,\ldots,D-1)$
\end{center}
is called the {\em standard sequence} corresponding to the eigenvalue $\theta$ (\cite[p.128]{bcn}).\\

Recall that a {\em clique} of a graph is a set of mutually adjacent vertices.  A clique $\mathcal{C}$ of a distance-regular graph with valency $k$, diameter $D\geq2$ and smallest eigenvalue $\theta_D$, is called {\em Delsarte clique} if $\mathcal{C}$ contains exactly $1-\frac{k}{\theta_D}$ vertices. It is known that for any clique $\mathcal{C}$ in $\Gamma$, the clique is Delsarte if and only if it is a completely regular code with covering radius $D-1$. A non-complete distance-regular graph $\Gamma$ is called {\em geometric} if there exists a set $C$ of Delsarte cliques such that each edge of $\Gamma$ lies in a unique $\mathcal{C}\in C$.\\

The {\em strong product} $G\boxtimes H$ of graphs $G$ and $H$ is a graph such that the vertex set of $G \boxtimes H$ is the Cartesian product $V(G)\times V(H)$ and any two different vertices $(u,v)$ and $(u',v')$ are adjacent in $G \boxtimes H$ if and only if ($u=u'$ or $u$ is adjacent to $u'$) and ($v=v'$ or $v$ is adjacent to $v'$). For a given positive integer $s$, the {\em $s$-clique extension} of a graph $G$ is the strong product $G\boxtimes K_s$ of $G$ and $K_s$, where $K_s$ is the complete graph (or clique) of size $s$.\\

A graph $\Gamma$ is said to be {\em of order $(s,t)$} if $\Gamma(x)$ is a disjoint union of $t+1$ cliques of size $s$ for every vertex $x$ in $\Gamma$. In this case, $\Gamma$ is a regular graph with valency $k=s(t+1)$ and every edge lies in a unique clique of size $s+1$.\\

The following lemma is straightforward.

\begin{lemma}{\rm (Cf. \cite[Proposition 4.3.2 and 4.3.3]{bcn} and \cite[Theorem 1.7.2]{gr})}\label{8}
Let $\Gamma$ be a distance-regular graph with valency $k$ and intersection number $c_2=1$. Then $\Gamma$ is of order $(s,t)$ and $b_1=t(a_1+1)=\frac{t}{t+1}k$, where $t\geq 1$ is an integer and $s=a_1+1$. Moreover if $t=1$ holds, then $\Gamma$ is a line graph.
\end{lemma}

A {\em Terwilliger graph} is a connected non-complete graph $\Gamma$ such that, for any two vertices $u, v$ at distance two, the subgraph induced on $\Gamma(u)\cap\Gamma(v)$ in $\Gamma$ is a clique of size $\mu$ (for some fixed $\mu \geq 1$). A {\em conference graph} is a strongly regular graph with parameters $v$, $k=\frac{v-1}{2}$, $a_1=\frac{v-5}{4}$ and $c_2=\frac{v-1}{4}$. An {\em antipodal graph}  is a connected graph $\Gamma$ with diameter $D>1$ for which being at distance 0 or $D$ is an equivalence relation. If, moreover, all equivalence classes have the same size $r$, then $\Gamma$ is also called an {\em antipodal $r$-cover}. \\

Let $x$ and $y$ be vertices of a distance-regular graph $\Gamma$. When $d(x,y)=t$ and $c_t=1$, we denote by $p[x,y]$ the unique shortest path connecting $x$ and $y$. Let $u,v$ be vertices of a distance-regular graph $\Gamma$ at distance $s$. Let $\Psi$ be the subgraph induced on $\Gamma_s(u)$ and let $\Psi_v$ be the connected component of $\Psi$ containing $v$. Let $\Delta(u,v)=\displaystyle\bigcup_{x\in\Psi_v}p[u,x]$.\\

For a vertex $x$ of a graph $\Gamma$, we write $\widehat{\Gamma(x)}$ for the set of vertices consisting of $x$ and its neighbors. Let $x\equiv y$ if $\widehat{\Gamma(x)}=\widehat{\Gamma(y)}$. Then, $\equiv$ is an equivalence relation, and we shall write $\overline{\Gamma}$ for the quotient $\Gamma/\equiv$ and $\overline{x}$ for the equivalence class of the vertex $x$. (I.e., $\overline{\Gamma}$ has vertices $\overline{x}$ for $x\in V(\Gamma)$ and $\overline{x}\sim_{\overline{\Gamma}}\overline{y}$ when $x \sim_{\Gamma} y$ and $\overline{x}\neq\overline{y}$.) $\overline{\Gamma}$ is called the {\em reduced graph} of $\Gamma$, and $\Gamma$ is called reduced when all equivalence classes have size one.\\

Recall the following interlacing result.

\begin{theorem}\label{1}{\rm (Cf.\cite{haem})}
Let $m \leq n$ be two positive integers. Let
$A$ be an $n\times n$ matrix, that is similar to a (real) symmetric matrix, and let
$B$ be a principal $m \times m$ submatrix of $A$. Then, for $i=1,\ldots , m$, $$\theta_{n-m+i}(A)\leq \theta_i(B)\leq \theta_i(A)$$
holds, where $A$ has eigenvalues $\theta_1(A) \geq \theta_2(A) \geq \ldots\geq \theta_n(A)$ and B has eigenvalues
$\theta_1(B) \geq \theta_2(B) \geq \ldots \geq \theta_m(B)$.\\
\end{theorem}

\section{A lower bound on the second largest eigenvalue}

In this section, we give a lower bound on the second largest eigenvalue and characterize the distance-regular graphs which attain this bound. This analysis will figure prominently in the proof of Theorem~\ref{6}.\\

Let $\Gamma$ be a distance-regular graph with valency $k$ and diameter $D$. Then the distinct eigenvalues of $\Gamma$ are those of the $(D+1)\times(D+1)$  tri-diagonal matrix\\ $$L_{\Gamma}=\left[
\begin{array}{cccccc}
 0 & k & & & &  \bf{0}\\
 c_1 & a_1 & b_1 & &  &\\
 & c_2 & a_2 & b_2 & &\\
 & & \ddots & \ddots & \ddots & \\
 & &   & c_{D-1} & a_{D-1}& b_{D-1} \\
 \bf{0}&&&& c_D & a_D
\end{array} \right],$$\\ where $c_i, a_i$ and $b_i$ $(0\leq i \leq D)$ are the intersection numbers of $\Gamma$.\\

 Let $L_{\Gamma}(i) :=\left[
\begin{array}{cccc}
 0 & k & &  \bf{0}\\
 c_1 & a_1 & b_1 &   \\
  & \ddots & \ddots & b_{i-1}  \\
  \bf{0}&& c_i & a_i
\end{array} \right]$ be the $(i+1)\times(i+1)$ upper left square
 submatrix of $L_{\Gamma}$, where $0\leq i\leq D$. Here note that the largest eigenvalue of $L_{\Gamma}(i)$, say $\mu_i$, is at least the average valency of the induced subgraph on $\{x\}\cup\Gamma(x)\cup\Gamma_2(x)\cup\cdots\cup\Gamma_i(x)$ for a vertex $x$ of the distance-regular graph $\Gamma$. i.e., $\mu_i>a_i+c_i$ holds.
\begin{lemma}\label{3}
Let $\Gamma$ be a distance-regular graph with diameter $D\geq 2t+2$ for some positive integer $t$. Then the second largest eigenvalue $\theta_1$ of $\Gamma$ is at least the largest eigenvalue $\mu_t$ of $L_{\Gamma}(t)$. Moreover, $\theta_1=\mu_t$ if and only if the graph $\Gamma$ is an antipodal distance-regular graph with diameter $D=2t+2$.
\end{lemma}

\begin{pf}
Let $x$ and $y$ be vertices of $\Gamma$ at distance $D(\geq2t+2)$. Then the induced subgraph of $\Gamma$ on $(\displaystyle\bigcup^t_{i=0}\Gamma_i(x)) \cup (\bigcup^t_{j=0}\Gamma_j(y))$ consists of two disjoint components and has the eigenvalue $\mu_t$ with multiplicity at least two. Thus, the inequality $\theta_1\geq\mu_t$ holds by Theorem~\ref{1}.

\ind Now we show that $\theta_1=\mu_t$ if and only if the graph $\Gamma$ is antipodal with $D=2t+2$. Let $(1, u_1, u_2, \ldots, u_D)$ be the standard sequence corresponding to $\theta_1$ and let $[1~~ u_1'~~ u_2'~~ \cdots~~ u_t']^T$ be the eigenvector of $L_{\Gamma}(t)$ corresponding to $\mu_t$. Then $u_1'=\frac{\mu_t}{k}$ and $c_iu_{i-1}'+a_iu_i'+b_iu_{i+1}'=\mu_tu_i'$ holds for $1\leq i \leq t-1$.\\

\ind If $\theta_1=\mu_t$, then clearly $u_i=u_i'$ for $1\leq i \leq t$ and hence $u_{t+1}=0$, as $c_tu_{t-1}'+a_tu_t'=\mu u_t'=\theta_1u_t=c_tu_{t-1}+a_tu_t+b_tu_{t+1}$ and $b_{t}>0$. Also we have $u_{t+1} = 0 > u_{t+2} > \ldots > u_D$ by \cite[p.130]{bcn}. The Perron-Frobenius Theorem \cite[Theorem 3.1.1]{bcn} allows us to conclude that $\theta_1$ is the largest eigenvalue of
$$\left[\begin{array}{cccc}
 a_{t+2} & b_{t+2} &  & {\bf 0}\\
 c_{t+3} & a_{t+3} & b_{t+3}  \\
 & \ddots & \ddots & b_{D-1} \\
{\bf 0} &  & c_D & a_D
\end{array} \right]$$
as the $b_i$'s and $c_i$'s are positive integers. Let $\Sigma(x)$ be the induced subgraph of $\Gamma$ on $\displaystyle\bigcup^D_{i=t+2}\Gamma_i(x)$ and let $z$ be a vertex of $\Gamma$ at distance $2t+2$ from $x$. Then the induced subgraph of $\Gamma$ on $\displaystyle\bigcup^t_{i=0}\Gamma_i(z)$, say $\Pi(z)$, is a subgraph of $\Sigma(x)$ with largest eigenvalue $\mu_t$. Let $\Lambda(z)$ be the connected component of $\Sigma(x)$ that contains $z$.  Then the following hold:\\

\ind (a) The largest eigenvalue of $\Lambda(z)$ is at most $\theta_1$;

\ind (b) $\Lambda(z)$ contains $\Pi(z)$ as an induced subgraph; and

\ind (c) $\Pi(z)$ has eigenvalue $\mu_t$. \\

This means that by the Perron-Frobenius Theorem,  $\Pi(z)$ must be equal to $\Lambda(z)$, and hence is a connected component of $\Sigma(x)$. This implies that $D=2t+2, c_D = k, b_{D-1} = 1$ and $a_1 = a_{D-1}$. As $\Pi(z)$ is connected, any vertex of $\Pi(z)\backslash\Gamma_t(z)$ has degree $k$ and any vertex of $\Gamma_t(z)$ has degree less than $k$. By \cite[Proposition 4.2.2 (ii)]{bcn}, we know that the graph $\Gamma$ is antipodal and that the converse is true. This completes the proof.
\end{pf}

 \textbf{Remark:} Lemma~\ref{3} was shown in \cite{kpy} for $t=1$.

\section{Terwilliger graphs}

The proof of Theorem~\ref{6} proceeds differently for different types of distance-regular graphs.  The type considered in this section is those with $c_2 > 1$ but which do not contain an induced quadrangle, commonly known as Terwilliger graphs.\

\begin{lemma}\label{4}
Let $m\geq2$ be an integer. Then there are only finitely many non-complete strongly regular Terwilliger graphs with smallest eigenvalue at least $-m$.
\end{lemma}

\begin{pf}
Let $\Gamma$ be a strongly regular graph and let $x$ and $y$ be vertices of $\Gamma$ at distance two. If $c_2=1$, then there is a maximal clique $C$ (with $a_1+2$ vertices) in $\Gamma_2(x)$ containing $y$. Then the distance $d(x,C)=2$. Hence, $C$ can not be a Delsarte clique, and $\Gamma$ is not geometric. So, we may assume $c_2\geq2$. If $\Gamma$ is geometric, then $\Gamma$ contains a quadrangle, as $c_2\geq2$ (\cite[Lemma 4.2]{kb}), and hence $\Gamma$ is not Terwilliger graph. By a result of Neumaier \cite{neu} (Cf.\cite[Theorem 1.1]{kb}), there are only finitely many non-geometric strongly regular graphs with smallest eigenvalue at least $-m$.
\end{pf}

As a consequence of Lemma~\ref{4}, we have

\begin{pro}\label{5}
Let $T\geq1$ be a real number. Then there are only finitely many distance-regular Terwilliger graphs with $c_2\geq2$ and second largest eigenvalue at least $\frac{b_1}{T}-1$.
\end{pro}

\begin{pf}
Let $\Gamma$ be a  distance-regular Terwilliger graph with $c_2\geq2$ and second largest eigenvalue at least $\frac{b_1}{T}-1$, and let $x$ be a vertex of $\Gamma$. Then the reduced graph $\overline{\Delta(x)}$ of $\Delta(x)$ is a strongly regular Terwilliger graph with smallest eigenvalue at least $-1-\frac{b_1}{\theta_1+1}\geq-1-T$, by \cite[Theorem 1.16.3]{bcn} and \cite[Theorem 4.4.3]{bcn}. By Lemma~\ref{4}, the number of vertices of $\Delta(x)$ is bounded. That is, the valency of $\Gamma$, say $k$, is bounded. As $c_2\geq2$, we have $h=1$, and hence the diameter $D$ of $\Gamma$ is at most $4^k$ by \cite[Theorem 5.9.8]{bcn}. This completes the proof.
\end{pf}

Note in \cite[Theorem 6.2]{kb}, they have a similar result when second largest eigenvalue is replaced by smallest eigenvalue.\\

In the case $T=2$, we can classify the corresponding Terwilliger graphs.\\

\begin{pro}\label{10}
Let $\Gamma$ be a distance-regular Terwilliger graph with $c_2\geq2$. If the second largest eigenvalue $\theta_1$ is bigger than $\frac{b_1}{2}-1$, then $\Gamma$
is either the icosahedron, the Conway-Smith graph, or the Doro graph.
\end{pro}

\begin{pf}
Let $x$ be a vertex of $\Gamma$. Then the local graph $\Delta(x)$ at $x$ is also a Terwilliger graph and has smallest eigenvalue at least
$-1-\frac{b_1}{\theta_1+1}>-3$ by \cite[Theorem 4.4.3]{bcn}. This implies that $\Delta(x)$ is an $s$-clique extension of a strongly regular Tewilliger graph
$\Sigma$ with parameters $(\overline{v},\overline{k},\overline{a_1},\overline{c_2})$ and smallest eigenvalue $\eta_{min}$ bigger than $-3$ by \cite[Theorem
1.16.3]{bcn} and Theorem~\ref{1}. If $\Sigma$ has a non-integral eigenvalue, then $\Sigma$ must be a conference graph with $\overline{c_2}=1$, as $\Sigma$ is a
Terwilliger graph. Thus, $\Sigma$ is the pentagon. As $\overline{a_1}=0$, we know $s=1$ by \cite[Theorem 1.16.3]{bcn}. The graph $\Gamma$ is therefore the icosahedron
by \cite[Proposition 1.1.4]{bcn}.

\ind From this point on, we will assume that $\Sigma$ has only integral eigenvalues. Then the smallest eigenvalue $\eta_{min}$ of $\Sigma$ is $-2$, as $\eta_{min}=-1$ implies
that $\Sigma$ is a complete graph. By \cite[Theorem 3.12.4]{bcn}, the graph $\Sigma$ is the Petersen graph. As $\overline{a_1}=0$, we know $s=1$ by \cite[Theorem
1.16.3]{bcn}. The graph $\Gamma$ is therefore either the Conway-Smith graph or the Doro graph by \cite[Theorem 1.16.5]{bcn}.
\end{pf}

\section{Distance-regular graphs with a fixed $\mathbf{\frac{b_t}{c_t}}$}

First, we will consider distance-regular graphs which have a fixed $\frac{b_t}{c_t}$ for some $t\geq1$. We will show that either the diameter or valency can be bounded in terms of $\frac{b_t}{c_t}$. We will isolate the special case of $\frac{b_2}{c_2}$, as in this case we may obtain an especially good diameter bound. These considerations will lead to the proof of Theorem~\ref{6}, given at the end of the section.

\begin{lemma}\label{7}
Let $C\geq\frac{1}{2}$ be a real number and $t$ be a positive integer. Let $\Gamma$ be a distance-regular graph with valency $k$ and diameter $D\geq t+1$. If  $\frac{b_t}{c_t}\leq C$, then either $D\leq 8C^2t$ or $k\leq2C$.
\end{lemma}

\begin{pf}
$\textbf{Claim 1}$: If $D\geq4t$, then $k\leq2C$ or $\frac{b_{4t}}{c_{4t}}\leq\frac{1}{2}C$.\\

\begin{description}
  \item[Case 1)] First consider the case $c_t\geq2$. Let $l$ be the minimum integer such that $c_l>c_{2t-1}$. Then by \cite[Theorem 1.1]{bhk}, we have $2t\leq l \leq (2t-1)+(2t-2)=4t-3$. As $c_l>c_{l-1}$ and $l\geq 2t$, we find $c_l\geq c_t+c_{l-t}\geq 2c_t$ by \cite[Proposition 1 (ii)]{koolen} and hence $\frac{b_l}{c_l}\leq\frac{1}{2}C$. So, $\frac{b_{4t}}{c_{4t}}\leq \frac{b_l}{c_l}\leq\frac{1}{2}C$.
\end{description}

From now on, we may assume $c_t=1$. i.e., $b_t\leq C$. In $\textbf{Case 2)}$ we consider the case $h\geq2$ and then we treat the case $h=1$ in $\textbf{Case 3)}$.

\begin{description}
  \item[Case 2)] Let us assume $h\geq2$. If $h\geq t$, then $b_1=b_t(\leq C)$ and $c_2=1$, and hence $k\leq2b_1\leq2C$ by Lemma~\ref{8}. We may therefore assume $t>h(\geq2)$.  As $h\geq2$,  by \cite[Theorem 2]{hira}, $c_{2h+3}\geq2$ and hence $\frac{b_{2h+3}}{c_{2h+3}}\leq\frac{1}{2}C$. Then $\frac{b_{4t}}{c_{4t}}\leq\frac{b_{2h+3}}{c_{2h+3}}\leq\frac{1}{2}C$, as $t>h(\geq2)$.
  \item[Case 3)] Now we assume $h=1$. If $c_{t+1}\geq2$, then $\frac{b_{4t}}{c_{4t}}\leq\frac{b_{t+1}}{c_{t+1}}\leq\frac{b_{t+1}}{2}\leq\frac{b_t}{2}\leq\frac{1}{2}C$, so we may assume $c_{t+1}=1$. We may also assume $c_{2t}=1$, as $c_{2t}\geq2$ implies $\frac{b_{4t}}{c_{4t}}\leq\frac{b_{2t}}{c_{2t}}\leq\frac{b_t}{2}\leq\frac{1}{2}C$. In \cite[Theorem 1.1]{hira2}, it was shown that if $c_{t+1}=1$ for some $t>0$, then for any two vertices $u$ and $v$ at distance $t$, there exists a distance regular graph $\Delta(u,v)$ with intersection array $\{a_t+1,a_t-a_1,\ldots,a_t-a_{t-1};1,1,\ldots,1\}$ as a subgraph of $\Gamma$, where the $a_i$'s are the intersection numbers of $\Gamma$.

      Let $x$ and $z$ be vertices of $\Gamma$ at distance $2t$. Then there is the unique vertex $y$ such that $d(x,y)=t$ and $d(y,z)=t$, as $c_{2t}=1=c_{2t-1}=\cdots=c_1$. This implies $V(\Delta(x,y))\cap V(\Delta(y,z))=\{y\}$. Let us consider the vertex $y$. Then the valency $a_t+1$ of the graph $\Delta(y,z)$ is at most the intersection number $b_t$ of $\Gamma$. This shows $a_t+c_t=a_t+1\leq b_t\leq C$. Thus, $k=a_t+c_t+b_t\leq 2b_t\leq 2C$.
\end{description}

$\textbf{Claim 2}$: If $k>2C$, then $D\leq8C^2t$.\\

Let $2^s\leq C <2^{s+1}$ for an integer $s\geq-1$ and assume $k>2C$. We will show by induction on $s$ that $D\leq 8C^2t$. If $s=-1$, then $\frac{1}{2}\leq C<1$, so $\frac{b_t}{c_t}\leq C<1$ and hence $D\leq t+t=2t$ by Lemma~\ref{pre} $(iii)$. Thus, \textbf{Claim 2} holds for $s=-1$. Now assume that \textbf{Claim 2} holds for $s$. Let $2^{s+1}\leq C<2^{s+2}$ and $k>2C$. If $D\leq 4t$, then $D\leq 8C^2t$ holds, as $C\geq 2^{s+1}=1$. We may therefore assume $D>4t$. Then, by \textbf{Claim 1}, $\frac{b_{4t}}{c_{4t}}\leq\frac{1}{2}C$. By the induction hypothesis, $D\leq8\times(\frac{1}{2}C)^2\times 4t=8C^2t$. This proves \textbf{Claim 2}, and \textbf{Claim 2} completes the proof of the lemma.
\end{pf}

The case $t=2, c_2 > 1$ may be of special interest, as in most cases one can improve the diameter bound markedly if desired.

\begin{lemma}\label{9}
Let $\Gamma$ be a distance-regular graph with diameter $D\geq3$ and $c_2\geq2$. If $\Gamma$ contains a quadrangle and $\frac{b_2}{c_2}<\frac{\alpha}{2}$ holds for some integer $\alpha\geq2$, then the diameter $D$ is at most $\alpha+1$.
\end{lemma}

\begin{pf}
First, we will consider the case $c_2\geq a_1$. As $\Gamma$ contains a quadrangle, the girth of $\Gamma$ is at most four, and this implies that the sequence $\{c_i\}_{i=1,\ldots,D}$ is a strictly increasing sequence by \cite[Theorem 5.2.5]{bcn}. Note that the inequality $c_i\geq\frac{i}{2}c_2$ holds for $i\geq2$ by \cite[Theorem 5.4.1]{bcn} and \cite[Proposition 1 (ii)]{koolen}. Thus, we have $\frac{b_2}{c_{\alpha}}\leq\frac{2}{\alpha}\times\frac{b_2}{c_2}<1$, and hence $D$ is at most $\alpha+1$ by Lemma~\ref{pre}$(iii)$.

\ind Now suppose $a_1>c_2$. As $\Gamma$ contains a quadrangle, we know that the inequality

\begin{center}
$ c_n-b_n\geq c_{n-1}-b_{n-1}+a_1+2\geq \cdots~~~~~~~~~~~~~~~~~~~~~~~~~$
\\ $~~~~~~~~~~~ \geq c_2-b_2+(n-2)a_1+(n-2)2>(n-1)c_2-b_2$
\end{center}

holds by \cite[Theorem 5.2.1]{bcn} for $n\geq3$. If $\frac{\alpha}{2}$ is an integer, then by putting $n-1=\frac{\alpha}{2}$, we have $c_{\alpha/2+1}>b_{\alpha/2+1}$, as $b_2<\frac{\alpha}{2}c_2$, which implies $D < \alpha+2$ by Lemma~\ref{pre}$(iii)$. Likewise, if $\frac{\alpha}{2}$ is not an integer, then we put $n-1=\frac{\alpha}{2}+\frac{1}{2}$. We then have $c_{n-1}-b_{n-1}>-\frac{c_2}{2}$ and $c_n-b_n>\frac{c_2}{2}$. Thus, $(c_{n-1}-b_n)+(c_n-b_{n-1})>0$. Therefore, at least one of $c_{n-1}-b_n$ and $c_n-b_{n-1}$ is positive, and $D$ can be at most $2n-2=\alpha+1$.
\end{pf}

We finally come to the proof of Theorem~\ref{6}.

\begin{pft2}
Let $C$ be a positive integer. Let $\Gamma$ be a distance-regular graph with valency $k$, diameter $D\geq 6$ and $\frac{b_2}{c_2}\leq C$. Then either $k\leq2C$ or $D\leq16C^2$ by Lemma~\ref{7}. There are only finitely many distance-regular graphs $\Gamma$ with $k\leq2C$ by \cite{bdkm}, so we may assume $D\leq16C^2$. Note that the second largest eigenvalue $\theta_1$ of $\Gamma$ is at least $\mu_2>a_2+c_2$ by Lemma~\ref{3}. Note further that the number $v$ of vertices of $\Gamma$ is at most $(3+C+C^2+\cdots+C^{D-2})k_2$, as $k_i=k\frac{b_1}{c_i}\frac{b_2}{c_2}\cdots\frac{b_{i-1}}{c_{i-1}}\leq k_2C^{i-2}$ for $2\leq i \leq D$ and $1+k\leq 2k_2$.

\ind If $c_2=1$, then $b_2\leq C$ and $\Gamma$ is a graph of order $(s,t)$. Thus, by Lemma~\ref{8}, either ($\Gamma$ is a line graph with $a_1=\frac{1}{2}k-1$) or $a_1\leq\frac{1}{3}k-1$. If $\Gamma$ is a line graph with $a_1=\frac{1}{2}k-1$, then by \cite[Theorem 4.2.16]{bcn}, $\Gamma$ is the flag graph of a regular generalized $2d$-gon of order $(s,s)$ for some $s>1$(which has diameter $2d$ and intersection array $\{2s,s,\ldots,s;1,1,\ldots,1\}$). Clearly in this case $s=b_2\leq C$ and $D=2d$. By \cite[Theorem 6.5.1]{bcn}, we know that $2d\in\{6,8\}$, since $s>1$, and hence there are only finitely many distance-regular graphs with $D\geq6$, $\frac{b_2}{c_2}\leq C$, $c_2=1$ and $a_1=\frac{1}{2}k-1$. We may therefore assume $a_1\leq\frac{1}{3}k-1$. As $a_2+c_2\leq C$ implies $k=a_2+b_2+c_2\leq 2C$, we may assume $a_2+c_2>C$. Hence, $\theta_1>\frac{1}{2}k$, as $b_2\leq C$ and $\theta_1>a_2+c_2>C$. Then the standard sequence ($u_i(\theta_1))_{i=0,\ldots,D}$ satisfies $u_1(\theta_1)=\frac{\theta_1}{k}>\frac{1}{2}$ and $b_1u_2(\theta_1)=((\theta_1-a_1)u_1(\theta_1)-1)$. So $b_1u_2(\theta_1)>\frac{1}{12}k-\frac{1}{2}\geq\frac{1}{24}k$, unless $k\leq12$, and thus $u_2(\theta_1)>\frac{1}{24}$, as $b_1<k$. Then the multiplicity $m_1$ of $\theta_1$ is at most  $\frac{v}{\sum^D_{i=0}u_i(\theta_1)^2k_i}<\frac{(3+C+C^2+\cdots+C^{D-2})k_2}{u_2(\theta_1)^2k_2}<\frac{(3+C+C^2+\cdots+C^{D-2})k_2}{(1/24)^2k_2}=576(3+C+C^2+\cdots+C^{D-2})$ by \cite[Theorem 4.1.4]{bcn}. As $k\leq \frac{(m_1+2)(m_1-1)}{2}$ (by \cite[Theorem 5.3.2]{bcn}), there are only finitely many distance-regular graphs with $D\geq6$, $\frac{b_2}{c_2}\leq C$, $c_2=1$ and $a_1\leq\frac{1}{3}k-1$.

\ind From this point forth, we assume $c_2\geq 2$. As $\frac{b_2}{c_2}\leq C$, it follows that $b_2\leq\frac{C}{C+1}k$ and hence $a_2+c_2\geq\frac{1}{C+1}k>\frac{1}{C+1}b_1-1$. First consider the case where $\Gamma$ is a Terwilliger graph. Let $T\leq C+1$ be the maximal integer such that $\theta_1>a_2+c_2\geq \frac{b_1}{T}-1$. Then there are only finitely many distance-regular Terwilliger graphs with $D\geq6$, $\frac{b_2}{c_2}\leq C$ and $c_2\geq2$ by Proposition~\ref{5}. We need therefore only consider the case where $\Gamma$ is not a Terwilliger graph, i.e., $\Gamma$ contains a quadrangle as an induced subgraph. As $D\geq 6$, we have $a_1<a_1+2\leq\frac{k+c_D}{6}\leq\frac{k}{3}$ by \cite[Corollary 5.2.2]{bcn}. If $b_2\leq \frac{1}{2}k$, then $\theta_1>a_2+c_2=k-b_2\geq\frac{1}{2}k$, and hence $\theta_1-a_1\geq \frac{1}{6}k$. In the same manner as above, we obtain $u_2(\theta_1)>\frac{1}{24}k$, unless $k\leq24$. Then $m_1\leq 576(3+C+C^2+\cdots+C^{D-2})$ and $k\leq \frac{(m_1+2)(m_1-1)}{2}$ hold, and therefore there are only finitely many distance-regular non-Terwilliger graphs with $D\geq6$, $\frac{b_2}{c_2}\leq C$, $c_2\geq2$ and $b_2\leq\frac{1}{2}k$. Thus, we may assume $b_2>\frac{1}{2}k$. Note that if we put $c_2+b_2=\alpha k$, for some $\frac{1}{2}<\alpha\leq1$, then $a_2=k-(b_2+c_2)=(1-\alpha)k$ and $\alpha k=c_2+b_2<(1+C)c_2$ hold. This implies $\theta_1>a_2+c_2>(1-\alpha)k+\frac{\alpha}{C+1} k\geq\frac{1}{C+1}k$.  As $c_2\geq\frac{b_2}{C}>\frac{k}{2C}$ and $b_1<k$, we have $k_2=k\frac{b_1}{c_2}<2Ck$. Hence we obtain $v\leq(3+C+C^2+\cdots+C^{D-2})k_2<2(3C+C^2+\cdots+C^{D-1})k$, which implies $m_1<\frac{v}{u_1(\theta_1)}\leq\frac{v}{1/(C+1)}\leq(C+1)(3+C+C^2+\cdots+C^{D-2})$. As $k\leq \frac{(m_1+2)(m_1-1)}{2}$, we see that there are only finitely many distance-regular non-Terwilliger graphs with $D\geq6$, $\frac{b_2}{c_2}\leq C$, $c_2\geq2$ and $b_2>\frac{1}{2}k$. This completes the proof of Theorem~\ref{6}.
\end{pft2}

We end this paper with the following question:\\

\textbf{Question:} Can Theorem~\ref{6} be strengthened as follows. Let $C$ and $t$ be positive integers. Are there only finitely many distance-regular graphs with valency at least three, diameter at least $2t+2$ and $\frac{b_t}{c_t}\leq C$?

\section{Acknowledgements}

 The second author was partially supported by the Basic Science Research Program through the National Research
Foundation of Korea(NRF) funded by the Ministry of Education, Science and Technology (Grant \# 2009-0089826). The third author was supported by Priority Research Centers Program through the National Research Foundation of Korea (NRF) funded by the Ministry of Education, Science and Technology (Grant \#2009-0094070).

\clearpage

\end{document}